\documentclass[12pt]{amsart}

\usepackage{amsmath, amscd, graphicx, latexsym, hyperref, times, rlepsf}

\oddsidemargin .25in \theoremstyle{plain}

\newtheorem{Thm}{Theorem}

\newtheorem{Lem}[Thm]{Lemma}

\newtheorem{Con}[Thm]{Conjecture}

\newtheorem{Cor}[Thm]{Corollary}

\newtheorem{Prop}[Thm]{Proposition}

\newtheorem{Rem}[Thm]{Remark}

\newcommand{\sn}{\operatorname{sn}}

\newcommand{\Hg}{\operatorname{Hg}}

\newcommand{\cg}{\operatorname{cg}}

\newcommand{\bfz}{{\mathbb{Z}}}

\newcommand{\OB}{\mathcal{OB}}

\def\v{\vskip.12in}

\begin{document}

\title[]{On the Heegaard genus of contact $3$-manifolds}

\author{Burak Ozbagci}

\address{Mathematical Sciences Research Institute \\ 17 Gauss Way \\ Berkeley, CA, 94720-5070}

\email{bozbagci@msri.org}

\subjclass[2000]{57R17}


\begin{abstract}

It is well-known that Heegaard genus is additive under connected sum of $3$-manifolds. We show that Heegaard genus of contact $3$-manifolds is not necessarily 
additive under \emph{contact} connected sum. 
We also prove some basic properties of the contact genus (a.k.a. open book genus \cite{rub}) of  $3$-manifolds, and compute this invariant for some $3$-manifolds.

\end{abstract}

\maketitle

\section{Introduction}

We assume that all $3$-manifolds are closed, connected and oriented and all contact structures are co-oriented and positive throughout this paper. 
Let $Y$ denote a $3$-manifold.  Given an open book $(B, \pi)$ on $Y$, where $B$ denotes the binding and $\pi$ denotes the fibration of  $Y-B $ over  $S^1$. 
It follows  that  $(\pi^{-1} ([0,1/2]) \cup B)$  and  $( \pi^{-1} ([1/2,1]) \cup B)$ are handlebodies which induce a Heegaard splitting of $Y$, where we view $S^1$ as the interval $[0,1]$ whose endpoints are identified with each other.  In this sense an open book is can be viewed as a special Heegaard splitting.   Note that a stabilization of an open book at hand corresponds to a stabilization of the induced Heegaard splitting. 

We define the Heegaard genus $\Hg(Y, \xi)$ of a contact $3$-manifold $(Y, \xi)$ as  the minimal genus of a Heegaard surface in any Heegaard 
splitting of $Y$ induced  from an open book supporting $\xi$. Equivalently, $\Hg(Y, \xi) =  1 + \sn(\xi) = \min \{ 1-\chi(\Sigma) \; | \;\Sigma $ is a page of an open book supporting $ \xi \}$, where $\sn(\xi)$  denotes the support norm of $\xi$ (cf. \cite{eto2}) and $\chi(\Sigma)$ denotes the Euler characteristic of $\Sigma$.  
This is certainly a generalization of the Heegaard genus adapted to contact 
$3$-manifolds.  It is well-known that Heegaard genus is additive under connected sum of $3$-manifolds. Here we show that 
Heegaard genus is sub-additive but not necessarily additive under connected sum of \emph{contact} $3$-manifolds.

Moreover we define the contact genus $\cg(Y)$ of a $3$-manifold $Y$ as the minimal Heegaard genus over all contact structures, i.e., $\cg(Y) = \min  \{  \Hg(Y, \xi) \;  |  \;\xi  $  is a contact structure on $Y\}$  
which, by Giroux's correspondence \cite{gi},  is the minimal genus of a Heegaard surface in any Heegaard 
splitting of $Y$ induced  from an open book. 
In other  words, the contact genus of a $3$-manifold is a topological invariant obtained by taking the minimum of the sum $2g+r-1$ over all open books, where $g$ and $r$ denote the genus of the page and the number of binding components  of the open book,  respectively.  We show that contact genus is sub-additive (and conjecture that it is additive) under  
connected sum of $3$-manifolds. 

We would like to point out that the contact  invariant was first studied by Rubinstein who named it the open book genus of $Y$ (cf. \cite{rub}).  We prefer to call it the contact genus to emphasize its connection with contact topology.  It is clear by definition that for any contact structure $\xi$ on $Y$ we have $$\Hg(Y) \leq \cg(Y) \leq \Hg(Y,\xi),$$ where $\Hg(Y)$ denotes the Heegaard genus of $Y$.  In \cite{rub}, it was shown that ``most'' $3$-manifolds of Heegaard genus $2$ have contact genus $>2$, which  implies the existence of $3$-manifolds where the first inequality above is strict.  In particular, it follows that not every Heegaard splitting of a $3$-manifold 
comes from an open book. 

Here we  show that  ``most'' $3$-manifolds of Heegaard genus $1$ have contact genus $>1$. Namely we show that a lens space which is not diffeomorphic to an oriented  circle bundle over $S^2$  have contact genus $ \geq 2 $. On the other hand, the contact genus of any oriented circle bundle over $S^2$  is equal its Heegaard genus.  We also show that there are many small Seifert fibered $3$-manifolds (which are not lens spaces) which have this property. Examples of such $3$-manifolds were constructed in \cite{rub}, but our examples are much simpler. We refer the reader to \cite {et2} and \cite {ozst} for more on open books and contact structures.

\section{Heegaard genus and contact connected sum}

Given any two contact $3$-manifolds $(Y_1, \xi_1)$ and $(Y_2, \xi_2)$. By removing a Darboux ball from each of these contact $3$-manifolds and gluing them along their convex boundaries by an orientation reversing map carrying respective characteristic foliations onto each other we get a well defined contact structure $ \xi_1 \# \xi_2 $ on the connected sum $Y_1 \#  Y_2$. The contact $3$-manifold 
$(Y_1 \# Y_2, \xi_1 \#  \xi_2)$ is called the contact connected sum of $(Y_1, \xi_1)$ and $(Y_2, \xi_2)$. It is well-known that Heegaard genus is additive under connected sum of smooth $3$-manifolds, which follows from Haken's Lemma. Here we show that

{\Thm\label{con} The Heegaard genus is sub-additive but not necessarily additive under connected sum of contact $3$-manifolds.}

\begin{proof} 
Let $\OB_i$ denote the open book 
realizing $\Hg(Y_i, \xi_i)$, for $i=1,2$.  Then the contact structure $ \xi_1 \# \xi_2$ on 
$Y_1 \#  Y_2$ is supported by the open book $\OB$ obtained by plumbing the pages of the open books $\OB_1$ and $\OB_2$ by Torisu \cite{to}.  
Denote a page of the open book $\OB_i$ by  $\Sigma_i$. It follows that $$ -\chi(\Sigma) =- \chi(\Sigma_1) - \chi(\Sigma_2)+1 ,$$
where $\Sigma $ denotes the page of the open book $ \OB$.  Thus we have $$ \Hg(Y_1 \# Y_2, \xi_1 \#  \xi_2) \leq \Hg(Y_1, \xi_1) + \Hg(Y_2, \xi_2), $$ which implies that  $\Hg$ is sub-additive under contact connected sum.

\v
Next we show  that $\Hg$ is not necessarily additive under contact connected sum. Let $\xi_d$ denote the overtwisted contact structure in $S^3$ whose $d_3$ invariant (cf. \cite{go}) is equal to the half integer $d$. The following result was obtained in \cite{dgs}: If  $(Y , \xi) $ is a contact structure with  $c_1 (\xi )$ torsion, then $$d_3 (Y , \xi \# \xi_d) = d_3 (Y , \xi) + d_3 (S^3 , \xi_d) + 1/2.$$  

Now suppose that $Y$ is an integral homology sphere. It follows that $c_1(\xi) =0$ for every 
contact structure $\xi$ on $Y$, and $Y$ carries a unique spin$^c$ structure. Thus  for an arbitrary contact structure $\xi$ on $Y$ we have $$d_3(Y, \xi \# \xi_{-\frac{1}{2}}) = d_3(Y, \xi) + d_3(S^3, \xi_{-\frac{1}{2}}) + \frac{1}{2} =  d_3(Y, \xi),$$ which implies 
that the connected sum $\xi \# \xi_{-\frac{1}{2}}$  is homotopic to $\xi$ as oriented plane fields (cf. \cite{go}). In fact,  $\xi \# \xi_{-\frac{1}{2}}$ is isotopic to $\xi$ by the classification of overtwisted contact structures due to Eliashberg \cite{eli}. As a consequence  we have 
$$ \Hg(Y, \xi \# \xi_{-\frac{1}{2}}) =  \Hg(Y, \xi).$$ 

On the other hand, in (\cite{eto2}, Lemma 5.5), it was proved that $\Hg(S^3, \xi_{-\frac{1}{2}}) =  2$.  Note that an open  
book realizing $\Hg(S^3, \xi_{-\frac{1}{2}})$ can be described  by taking a pair of pants as a page  and $t_1 t_2^{-2} t_3^{-3}$ as the monodromy,  where $t_i$ denotes a right-handed Dehn twist along a boundary component. Consequently we have

$$ \Hg(Y \# S^3, \xi \# \xi_{-\frac{1}{2}}) < \Hg(Y, \xi) +  \Hg(S^3, \xi_{-\frac{1}{2}}). $$

 \end{proof}

 



\section{Contact genus of three dimensional manifolds}

Here we provide some basic properties of the contact genus of  $3$-manifolds,  and compute this invariant for some $3$-manifolds.

{\Prop \label{basic}    Let $Y$ denote  a $3$-manifold. Then we have

\begin{itemize}

\item[($a$)] $\cg(Y)  \geq 0$ ($ = 0 $ if and only if  $Y \cong S^3$),

\item[($b$)] $  \cg(Y)  = 1 $ if and only if $Y$ is an oriented  circle bundle over $S^2$ (which is not diffeomorphic to $S^3$).

\end{itemize}

} 

\begin{proof}  For a $3$-manifold $Y$, $\cg(Y)$ is obtained by taking the minimum of the sum $2g+r-1$ over all open books, where $g$ and $r$ denote the genus of the page and the number of binding components  of an open book,  respectively. Hence  we have $0 \leq \cg(Y) $ for an arbitrary $3$-manifold $Y$, since $g \geq 0$ and $r \geq 1$.  It is clear that the absolute minimum of the expression $2g+r-1$ is realized when 
$g=0$ and $r=1$ and  the open book with disk pages and trivial monodromy supports the unique tight contact structure on $S^3$, which proves ($a$). 

To prove ($b$),  we note that $\cg(Y) = 1 $ is realized if and only if $g=0$ and $r=2$. Any 
self-diffeomorphism of  an annulus is given by $t_c^m$, for some $m \in \bfz$, where $c$ is the core of the annulus, and $t_c$ denotes a right-handed Dehn twist along $c$.  If $m \geq 0$, this open book supports the unique tight contact structure on the lens space $L(m, -1)$ which is an oriented circle bundle over $S^2$ with Euler number $m$. Otherwise (i.e., when $m <0$) the induced contact structure is the overtwisted contact structure on $L(-m, 1)$ which is an oriented circle bundle over $S^2$ with Euler number $m$.  Combining, we showed that  $  \cg(Y)  = 1 $ if and only if $Y$ is an oriented  circle bundle over $S^2$, which is not diffeomorphic to $S^3$.   

\end{proof}

Note that oriented circle bundles over $S^2$ are very special lens spaces and therefore we immediately conclude from Proposition~\ref{basic}  that

{\Cor \label{str}  Most  $3$-manifolds of Heegaard genus $1$ have contact genus $>1$. }

\v  
 
For example, $\cg(L(5,3))=2$, since $L(5,3)$ is not a  circle bundle over $S^2$ and it carries a (tight)  
contact structure which is supported by a planar  open book with three binding components.

{\Lem \label{ba}  We have  $ \cg(Y_{p,q,r})  \leq  2$ , where $Y_{p,q,r}$ denotes the $3$-manifold depicted in Figure~\ref{seif}, with $p,q,r \in \bfz$. Moreover if $|p| >1$, $|q| > 1$ and $|r| > 1$ then $ \cg(Y_{p,q,r}) =2 $. }

\begin{proof}
 
It follows from \cite{eto2} that $Y_{p,q,r}$ has a planar open book with at most three binding components, which indeed proves that $ \cg(Y_{p,q,r})  \leq 2 $. Moreover, under the assumption that $|p| >1, |q|>1$, and $|r| >1$, the $3$-manifold $Y_{p,q,r}$ is not diffeomorphic to any lens space and hence $ \cg(Y_{p,q,r})=2 $ by Proposition~\ref{basic}.  

\begin{figure}[ht]
  \relabelbox \small {\epsfxsize=2.5in
  \centerline{\epsfbox{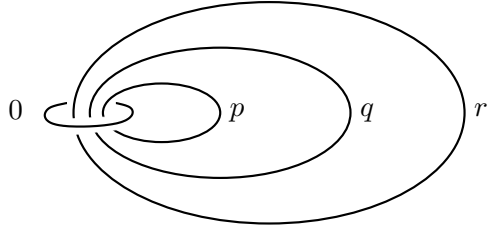}}}

\relabel{p}{$p$}

\relabel{q}{$q$}

\relabel{r}{$r$}

\relabel{s}{$0$}

\endrelabelbox
        \caption{Integral surgery diagram for the small Seifert fibered $3$-manifold $Y_{p,q,r}$} \label{seif}
 \end{figure}

\end{proof}

When we drop the assumption on $p,q$ and $r$ in Lemma~\ref{ba}, we observe that $Y_{p,q,r}$ is diffeomorphic to either $S^3$, $S^1 \times S^2$, a lens space, or certain connected sums of these for some values of the integers $p,q$ and $r$.

{\Rem Note that Lemma~\ref{ba} exhibits examples of $3$-manifolds $Y= Y_{p,q,r}$ for which $\Hg(Y)=\cg(Y)=2$, although most $3$-manifolds of Heegaard genus $2$ have contact genus $> 2$ as was shown by Rubinstein \cite{rub}.}

{\Lem \label{es} We have $ \cg (\#_k S^1 \times S^2) = k$, for $k \geq 1$.}

\begin{proof}

Since $\Hg (\#_k S^1 \times S^2) = k$, we know that $\cg (\#_k S^1 \times S^2) \geq k$. Hence to show that $\cg (\#_k S^1 \times S^2) = k$, we just need to realize this lower bound by a Heegaard splitting of $\#_k S^1 \times S^2$ induced from an open book. We use the fact that the unique tight contact structure on $\#_k S^1 \times S^2$ is supported by an planar open book with $k+1$ binding components, whose monodromy is the identity map. 

\end{proof}

 The proof of the following result is similar to the proof of Theorem~\ref{con}.

{\Prop  Let $Y_i$ denote  a $3$-manifold, for $i=1,2$. Then we have $$ \cg(Y_1 \#  Y_2) \leq \cg(Y_1) + \cg(Y_2) .$$}

{\Con Contact genus is additive under connected sum of $3$-manifolds.}

\v  Note that if $\Hg(Y_i) = \cg(Y_i)$ for $i=1,2$, then $\cg(Y_1 \# Y_2) = \cg(Y_1) + \cg(Y_2)$.

\v \noindent {\bf {Acknowledgement}}: The author would like to thank John B. Etnyre and Ko Honda for helpful conversations and the Mathematical Sciences Research Institute for its hospitality during the \textit{Symplectic and Contact Geometry and Topology} program 2009/10. B.O. was partially supported the Marie Curie International Outgoing Fellowship 236639.

\end{document}